\documentclass[letterpaper,journal,twocolumn,twoside,final,9pt]{IEEEtran}
\usepackage{amssymb,amsfonts}

\providecommand\href[2]{#2}
\newtheorem {theorem}{Theorem}
\newtheorem {corollary}{Corollary}

\newtheorem {proposition}{Proposition}
\newtheorem {remark}{Remark}
\newtheorem {constr}{Construction}
\newtheorem {conjecture}{Conjecture}

\newcommand\F[1][n]{F^{#1}_{\scriptscriptstyle\rm E\hspace{-0.1ex}v}}
\newcommand\FF[1][n]{F^{#1}_{\scriptscriptstyle\rm O\hspace{-0.1ex}d}}

\def\eqdf{\triangleq}

\begin{document}
\def\hrefDAOR#1{\href{http://math.nsc.ru/conference/DAOR'04/main.html}{#1}}
\title{On the Number of $1$-Perfect Binary Codes:\\ a Lower Bound%
\thanks{
This is author's version of the correspondence in the
\href{http://ieeexplore.ieee.org/xpl/RecentIssue.jsp?punumber=18}{IEEE Transactions on Information Theory}
54(4) 2008, 1760-1765,
Digital Object Identifier
\href{http://dx.doi.org/10.1109/TIT.2008.917692}{10.1109/TIT.2008.917692},
\copyright 2008 IEEE.
Personal use of the material of the correspondence is permitted.
However, permission to reprint/republish
the material for advertising or promotional
purposes or for creating new collective works
for resale or redistribution to servers or lists,
or to reuse any copyrighted component of this work
in other works must be obtained from the IEEE.}
\thanks{The material in this correspondence was presented in part at the
10th International Workshop on Algebraic and Combinatorial Coding Theory
ACCT-10, Zvenigorod, Russia, September 2006 
and, in a brief form,
at the Russian conference
``\hrefDAOR{{D}iscrete Analysis and} \hrefDAOR{Operations Research}''
\hrefDAOR{{DAOR}'2004},
Novosibirsk, Russia, June-July 2004, \href{http://www.math.nsc.ru/conference/DAOR'04/daor4.pdf}{p95}.}
}
\author{\href{http://arxiv.org/find/grp_math/1/au:+Krotov_Denis/0/1/0/all/0/1}{Denis~S.~Krotov} and
\href{http://arxiv.org/find/grp_math/1/au:+Avgustinovich_Sergey/0/1/0/all/0/1}{Sergey~V.~Avgustinovich}
\thanks{The authors are with the Sobolev Institute of Mathematics,
pr-t Ak. Koptyuga 4, Novosibirsk, 630090, Russia
(e-mail: krotov@math.nsc.ru, avgust@math.nsc.ru)}
 }
\maketitle
\begin{abstract}
We present a construction of $1$-perfect binary codes,
which gives a new lower bound on the number of such codes.
We conjecture that this lower bound is asymptotically tight.

\end{abstract}

\newcommand{\EExpLog}[2][0ex]{{{\raisebox{#1}{$2$}}^{\frac{n}{#2}-{\log}\frac{n}{#2}-1}}}
\newcommand{\EExp}[2][0ex]{{{\raisebox{#1}{$2$}}^{\frac{n}{#2}-1}}}
\newcommand{\EexpLog}[2][-0.2ex]{{{\raisebox{#1}{$\scriptstyle 2$}}^{\frac{n}{#2}-{\log}\frac{n}{#2}-1}}}
\newcommand{\Eexp}[2][-0.2ex]{{{\raisebox{#1}{$\scriptstyle 2$}}^{\frac{n}{#2}-1}}}
\newcommand{\ExpHam}[2][0ex]{2^{\EexpLog[-0.2ex]{#2}}}

\section{ Introduction}
The paper is devoted to the problem of enumeration of $1$-perfect binary codes.
Such codes, as any optimal codes,
are extremal objects of the theory of error-correcting codes.
In addition, perfect codes are a special type of combinatorial configurations.
The construction of $1$-perfect binary codes presented in the current paper gives
the most powerful known class of such codes and leads to a lower bound on their number.

The first known construction \cite{Vas:nongroup_perfect}
of nonlinear $1$-perfect binary codes gives the lower bound
\begin{eqnarray} \nonumber
B(n-1)&\geq&\ExpHam{2}\cdot\ExpHam{4}\cdot\ExpHam{8}\cdot      \dots
\end{eqnarray}
on the number $B(n-1)$ of $1$-perfect codes of length $n-1=2^m-1$ (here
and in what follows ${\log}$ means ${\log}_2$).
This bound was improved in \cite{AvgSol:alpha-comp} and \cite{Mal:lb},
where some useful ideas exploited 
in this paper were proposed.
The best known lower bound \cite{Kro:quas_l_b} is
\begin{eqnarray} \label{eq:0}
B(n-1)&\geq&\ExpHam{2}\cdot 3^{\Eexp{4}}\cdot \ExpHam{4}.
\end{eqnarray}
The result of \cite{Kro:quas_l_b} was formulated in the terms of a partial case \cite{Phelps84}
of the generalized concatenation construction (see, e.g., \cite{ZinLob:2000}),
which allows to construct $1$-perfect binary codes from distance $2$
$q$-ary MDS codes (only the case $q=4$ is useful for the lower bound),
or $n$-ary quasigroups (of order $4$).
The lower bound on the number of
$n$-ary quasigroups of order $4$ given in \cite{Kro:quas_l_b}
is asymptotically tight \cite{KroPot:Lyap,PotKro:asymp};
therefore, such a way to evaluate the number
of $1$-perfect binary codes has been exhausted.

The best known upper bound \cite{Avg:ub} on the number of $1$-perfect binary codes
is
$
2^{2^{n-(3/2){\log} n+\log{\log} (e n)}}.
$

The local-automorphism method presented in this work is a further development
of the methods \cite{Vas:nongroup_perfect,AvgSol:alpha-comp,Mal:lb,Kro:quas_l_b}.
Since there is one-to-one correspondence between $1$-perfect binary codes
and extended $1$-perfect binary codes, the results are formulated in terms of
extended $1$-perfect binary codes.

\section{ Preliminaries}
Let $F^n$ ($\F$, $\FF$) be the set of the binary $n$-words
(with even or odd number of ones, respectively) with the Hamming distance $d$
and mod~$2$ coordinate-wise addition. Given $\bar x \in F^n$, put $wt(\bar x) \eqdf d(\bar x,\bar 0)$.
The {\it neighborhood} of $S\subseteq F^n$ is the set
$\Omega(S)\eqdf \bigcup_{\bar x\in S}\Omega(\bar x)$ where
$\Omega(\bar x)\eqdf \{\bar y\in F^n \,|\,  d(\bar y,\bar x)=1\}$.

A set $C\subseteq F^n$ is called a \emph{distance $d$ code} (\emph{of length $n$})
{{if}} the Hamming distance between any two different words in $C$ is not less than $d$.
An {\it extended $1$-perfect code\,} is a set $C\subseteq \F$
such that the neighborhoods of the words of $C$
are pairwise disjoint and $\Omega(C)=\FF$.
It follows that $C$ is a distance $4$ code
of cardinality $|C|=|\FF|/n=2^{n-{\log} n-1}$ and $n$ is a power of $2$.
On the other hand, if the code $C\subset \F$ has distance $4$ and $|C|=2^{n-{\log} n-1}$, then, obviously, $C$ is an extended $1$-perfect code.
In what follows we assume $n=2^m\geq 16$.

The following formulas define some useful sets $V^t,A^t\subset\F$ and
give a representation of the extended Hamming code:
\begin{eqnarray}
V^t & \eqdf & \{(\bar{v},\bar v,0,\ldots ,0)\in F^n\,|\,\bar v\in \F[2^{m-t}]\},
\\ \nonumber&&
 t=1,\ldots,m-1
\\
 A^1 & \eqdf & V^1
\nonumber \\
A^{t} & \eqdf & V^t + A^{t-1} = \bigcup_{\bar r\in V^{t}}(\bar r+A^{t-1}), \label{eq:At}
\\ \nonumber&&
 t=2,\ldots,m-1
 \\
H & \eqdf & A^{m-1}. \label{eq:1}
\end{eqnarray}
The following is straightforward:
\begin{proposition}
The set $H$ defined by {\rm (\ref{eq:1})} is a linear extended $1$-perfect code,
i.\,e., the extended Hamming code (see e.\,g. {\rm\cite[\S1.7]{MWS}}),
which is the only linear extended $1$-perfect code,
up to coordinate permutation.
\end{proposition}

We say that a set $G\subset \F$ is a \emph{component of order $t\in\{1,\ldots,m-1\}$}
{{if}} $|G|=|A^t|$ and $\Omega(G)=\Omega(A^{t})$.
Let $t\in\{1,\ldots,m-1\}$ and $\bar \mu\in F^{2^{m-t}}\times {0}^{n-2^{m-t}}$;
we say that a set $M\subset \F$
is a \emph{$\bar \mu$-component} (of order $t$)
{{if}} $M=G+\bar \mu$ for some component $G$ of order $t$.

Let ${\rm Aut}(F^n)$ denote the group of isometries $F^n$
(it coincides with the automorphism group of the distance-one graph of $F^n$).
It is known that each isometry $g\in {\rm Aut}(F^n)$ has a unique
representation $g(\cdot) = \bar v+\pi(\cdot)$ where $\bar v$ is a \emph{shift vector}
 from $F^n$ and $\pi$ is a coordinate permutation.
If $\pi=Id$, i.\,e., $g(\cdot) = \bar v+\cdot$, then the isometry $g$ is called a \emph{translation}.
If $S\subseteq F^n$, then ${\rm Aut}({S})$ is the group of isometries $g$ of $F^n$ such that $g(S)=S$.
For a collection ${\mathbf{S}}=\{S_1,\ldots ,S_l\}$ of subsets of $F^n$,
by ${\rm Aut}(\mathbf{S})$ denote the group of isometries $g$ of $F^n$ such that
for each $S\in\mathbf{S}$ the set $g(S)$ is also in $\mathbf{S}$.
In what follows we use calligraphic letters to denote subgroups of ${\rm Aut}(F^n)$.
Put $${\cal A}^t \eqdf {\rm Aut}(\Omega(A^{t}))$$
where $A^t$ is specified by (\ref{eq:At}).

\begin{proposition}\label{pro:1} {\rm (a)} Let $t\in\{2,\ldots,m-1\}$ and
let for each $\bar\mu\in V^t$ the set
$B_{\bar\mu}$ be a ${\bar\mu}$-component of order $t-1$.
Then the set $B=\bigcup_{\bar\mu\in V^t} B_{\bar\mu}$ is a component of order $t$.

{\rm (b)} If $B$ is a component of order $t$ and $g^t\in {\cal A}^t$,
then $g^t(B)$ is a component of order $t$ too.

{\rm (c)} A component of order $m-1$ is an extended $1$-perfect code.
\end{proposition}
\begin{proof}
(a) follows from (\ref{eq:At}), the definition of a component of order $t$,
and the equality $\Omega(\bar\mu+A^{t-1})=\Omega(B_{\bar\mu})$,
which holds by the definition of $\bar\mu$-component of order $t-1$.

(b) and (c) are straightforward from the definitions.
\end{proof}

\section{ The LA Construction of Extended $1$-Perfect Codes}
Proposition \ref{pro:1} is all we need to see that the following construction
leads to an extended $1$-perfect code.
The idea of the construction is to take the Hamming code and apply
isometries of $F^n$ to parts of the code in such a way that the neighborhood of the part
does not change. We call such isometries \emph{local automorphisms}; a local automorphism
acts on a part of the code and does not change the neighborhood of that part.
At the first stage we take components of order $1$ as such parts;
at the second stage, components of order $2$; and so on.
At the last stage, we ``turn'' the whole code.

\begin{constr}[LA -- local automorphisms] \label{con:1}
Assume for any integer $t \in \{2,\ldots, m\}$ and for any words
${\bar r_{i}\in V^{i}}$, $i=t,\ldots,m-1$
we have a local automorphism $g_{\bar r_t,\ldots ,\bar r_{m-1}}\in {\cal A}^{t-1}$;
in particular, $g\in {\cal A}^{m-1}$.
Then (as follows by induction on $t$ from Proposition \ref{pro:1})
the set $C$ represented by the following formulas is an extended $1$-perfect code.
\begin{IEEEeqnarray}{rCl}
A^1_{\bar r_{2},\ldots ,\bar r_{m-1}} & \eqdf & V^1,\qquad \bar r_{i}\in V^i
\nonumber\\
A^{t}_{\bar r_{t+1},\ldots ,\bar r_{m-1}} & \eqdf & \bigcup_{\bar r_{t}\in V^{t}}
\left(\bar r_t+g_{\bar r_t,\ldots ,\bar r_{m-1}}(A^{t-1}_{\bar r_t,\ldots ,\bar r_{m-1}})\right),
\nonumber \\ \nonumber&&
 \qquad\qquad t=2,\ldots,m-1
\label{eq:2}
\\
C & \eqdf & g(A^{m-1}).
\end{IEEEeqnarray}
\end{constr}

In Construction \ref{con:1} each code can be obtained in more than one way.
To evaluate the number of the codes that can be constructed in this way we need
 stronger restrictions on local automorphisms.

Let
\begin{IEEEeqnarray*}{rCl}
{\cal B}^1 &\eqdf& {\rm Aut}(A^{1})
\\
{\cal B}^t &\eqdf& {\rm Aut}\big(\{\bar r+\Omega(A^{t-1})\}_{\bar r\in V^{t}}\big),\qquad t=2,\ldots ,m-1.
\end{IEEEeqnarray*}
For each $t=1,\ldots ,m-1$ we fix a set ${\cal D}^{t}$ of representatives of
cosets from ${\cal A}^{t}/{\cal B}^{t}$.
Moreover, we choose the representatives in such a way that the following holds:
for two cosets $D_1,D_2\in{\cal A}^{t}/{\cal B}^{t}$ and their representatives
$d_1,d_2\in {\cal D}^{t}$, $d_1 \in D_1$, $d_2\in D_2$, the equality
$D_1=\tau D_2$ with some translation $\tau$ implies $d_1=\tau d_2$
(this condition is essential for the definition of degenerate collection
and Proposition~\ref{p:bold} below).

It can be shown by induction that

\begin{proposition}
The restrictions
$
g_{\bar r_t,\ldots ,\bar r_{m-1}}\in {\cal D}^{t-1}
$
do not reduce the set of codes that can be represented by {\rm (\ref{eq:2})}.
\end{proposition}
\begin{proof}
Assume that
$$G^t \eqdf \bigcup_{\bar r \in V^t}\left(\bar r+g_{\bar r}(G^{t-1}_{\bar r})\right)$$
where $G^{t-1}_{\bar r}$ is a component of order $t-1$
and $g_{\bar r}\in{\cal A}^{t-1}$
for all $\bar r\in V^t$.
Let $g\in {\cal A}^t$ and $g=dh$ where $d\in{\cal D}^t$
and $h\in {\cal B}^t$.

We claim that
\begin{equation}\label{eq:cl}
    g (G^t) = d (G'^t) \quad \mbox{where} \quad G'^t \eqdf \bigcup_{\bar q \in V^t}\left(\bar q+g'_{\bar q}(G^{t-1}_{\rho\bar q})\right)
\end{equation}
for some $g'_{\bar q}\in{\cal A}^{t-1}$ and permutation $\rho:V^t\to V^t$.
Indeed, by the definition of ${\cal B}^t$, for all $\bar r\in V^t$ we have
$$
h(\bar r+\Omega(A^{t-1}))=\rho^{-1}\bar r+\Omega(A^{t-1})
$$
where $\rho$ is some permutation on $V^t$.
So, we see that $h_{\bar r}(\cdot)\eqdf \rho^{-1}\bar r+h(\bar r+\cdot)$
belongs to ${\cal A}^{t-1}$.
Then, replacing $\bar r$ by $\bar q\eqdf \rho^{-1}\bar r$,
we see that (\ref{eq:cl}) holds with
$g'_{\bar q} \eqdf h_{\rho\bar q} g_{\rho\bar q}$.

So, using (\ref{eq:cl}), we can step-by-step replace the operators
$g_{...}\in {\cal A}^t$ by  $d_{...}\in {\cal D}^t$,
starting from $t=m-1$ and finishing with $t=1$.
\end{proof}

Therefore the following construction gives the same set of codes
as Construction \ref{con:1}.
\begin{constr}[LA, upper bound] \label{con:2}
Assume that for any integer $t \in \{2,\ldots, m\}$ and for any words
${\bar r_{i}\in V^{i}}$, $i=t,\ldots,m-1$,
we have $g_{\bar r_t,\ldots ,\bar r_{m-1}}\in {\cal D}^{t-1}$;
in particular, $g\in {\cal D}^{m-1}$.
Then the set $C$ defined by the formulas {\rm (\ref{eq:2})} is an extended $1$-perfect code.
\end{constr}

As we will see below (Theorem~\ref{th:Kla}),
almost all ($n\to\infty$) codes represented by Construction \ref{con:2}
have a unique representation and this gives a good upper estimation
\begin{equation}\label{up}K_{\scriptscriptstyle LA}(n)\leq
|{\cal D}^{m-1}|\prod_{t=1}^{m-2}
|{\cal D}^t|^{|V_{t+1}|\cdot|V_{t+2}|\cdot\ldots \cdot|V_{m-1}|}
\end{equation}
for the number $K_{\scriptscriptstyle LA}(n)$ of different extended
$1$-perfect codes of length $n$
obtained by the method of local automorphisms (\emph{LA}), i.\,e., by
Construction \ref{con:1} or \ref{con:2}.
To show that the number of different LA codes is close to this value,
we need some more restrictions on $g_{\bar r_t,\ldots ,\bar r_{m-1}}$.

Assume $L$ is a linear subspace of $F^n$ and for each $\bar r\in L$ we have
$g_{\bar r}\in {\cal D}^{t}$ and
$g_{\bar r}(\cdot)= \bar v_{\bar r}+\pi_{\bar r}(\cdot)$.
We say that the collection $\{g_{\bar r}\}_{\bar r\in L}$ is \emph{degenerate}
{{if}} the following conditions hold:

$\bullet$ the permutation $\pi_{\bar r}$ does not depend on $\bar r$,
i.\,e., $\pi_{\bar r}=\pi$ for all $\bar r\in L$;

$\bullet$ the set $\{\bar r + \bar v_{\bar r} \,|\, {\bar r\in L} \}$ is an affine
subspace of $F^n$.

Otherwise we say that $\{g_{\bar r}\}_{\bar r\in L}$ is \emph{nondegenerate}.

\begin{constr}[LA, lower bound] \label{con:3}
In addition to the conditions of Construction \ref{con:2} we require that
the collection $\mathbf{g}_{\bar r_{t+1},\ldots ,\bar r_{m-1}}=\{g_{\bar r_t,\ldots ,\bar r_{m-1}}\in {\cal D}^{t-1}\}_{\bar r_t\in V^t}$
is nondegenerate for every
$t\in \{2,\ldots,{m-1}\}$, $\bar r_{t+1}\in  V^{t+1}$, \ldots, $\bar r_{m-1}\in  V^{m-1}$.
\end{constr}

\section{ Calculations}

In this section we establish some facts concerning
the structure of order-$t$ components and related objects,
on which the main result is based.
Given $G\subseteq \F$, put
$$\Theta(G) \eqdf \{\bar x\in \F\,|\, \Omega(\bar x)\subseteq \Omega(G)\};$$
clearly, $G\subseteq\Theta(G)$ and $\Omega(\Theta(G))=\Omega(G)$.
The following fact is also straightforward:
\begin{proposition} \label{p:Theta} For each $G,G'\subseteq \F$ the equality
$\Omega(G)=\Omega(G')$ means $\Theta(G)=\Theta(G')$ and vice~versa.
\end{proposition}
For each $t=1,\ldots,m$ and $\bar x=(\bar x_0,\ldots,\bar x_{2^t-1})\in (F^{2^{m-t}})^{2^t}=F^n$
define the generalized parity check
$$p^t(\bar x) \eqdf \sum_{i=0}^{2^t-1}\bar x_i.$$

\begin{proposition} \label{p:At}
Let $1\leq t \leq m-1$; then the following claims hold:\\
{\rm (a)} $p^t(\bar x)=\bar 0$ for all $\bar x\in A^t;$\\
{\rm (a')} $|A^t|=2^{2^{m-t}(2^t-1)-t};$\\
{\rm (b)} $\Omega(A^t)=\{\bar x\in F^n\,|\,wt(p^t(\bar x))=1\};$\\
{\rm (b')} $|\Omega(A^t)|=2^{2^{m-t}(2^t-1)+m-t};$\\
{\rm (c)} if $t<m-1$, then $\Theta(A^t)=\{\bar x\in F^n\,|\,p^t(\bar x)=\bar 0\};$\\
{\rm (c')} if $t<m-1$, then $|\Theta(A^t)|=2^{2^{m-t}(2^t-1)};$\\
{\rm (c'')} $\Theta(A^{m-1})=\F$.
\end{proposition}
\begin{proof}
(a) and (a') are straightforward from the definition of $A^t$.

(b') Since the code distance of $A^t$ is $4$, we have $|\Omega(A^t)|=n|A^t|$.

(b) It follows from (a) that $wt(p^t(\bar x))=1$ for all $\bar x\in\Omega(A^t)$.
On the other hand, by (b'), we have $|\Omega(A^t)|=|\{\bar x\in F^n\,|\,wt(p^t(\bar x))=1\}|$.

(c) It follows from (b) that $p^t(\bar x)=\bar 0$ implies $\bar x\in \Theta(A^t)$.
Assume $p^t(\bar x)\neq\bar 0$. If $t<m-1$, then there is
$\bar y\in \Omega(\bar x)$ such that $wt(p^t(\bar y))>1$; therefore, $\bar x\not\in\Theta(A^t)$.

(c') follows from (c).

(c'') It follows from (b) or (b') that $\Omega(A^{m-1})=\FF$; thus $\Theta(A^{m-1})=\F$.
\end{proof}

In what follows we will use the `array' representation of elements of $F^n$:
$$\bar x=(x^t_{0,0},...,x^t_{0,2^{m-t}-1},x^t_{1,0},
\ldots,
x^t_{2^{t}-1,2^{m-t}-1})=
(x^t_{i,j})_{i,j}$$
where indexes $i,j$ change in lexicographical order.
I.\,e., for each $t=1,\ldots,m-1$ an element $\bar x$ in $F^n$
can be viewed as $2^{t}\times 2^{m-t}$-array
$$
\left(
\begin{array}{l@{\ }l@{\ }l@{\ \ }l}
x^t_{0,0} & x^t_{0,1} & \ldots & x^t_{0,2^{m-t}-1}\\
\ldots &\ldots &\ldots &\ldots \\
x^t_{2^{t}-1,0}& x^t_{2^{t}-1,1}& \ldots & x^t_{2^{t}-1,2^{m-t}-1}
\end{array}
\right)
$$
In these terms, $p^t(\bar x)$ is the  sum of rows of $(x^t_{i,j})_{i,j}$.
For further calculations, we introduce the sets
\begin{IEEEeqnarray*}{rCl}
  {B}^1 & \eqdf & V^1  \\
  {B}^t & \eqdf & V^t + \Theta(A^{t-1}) 
  , \qquad t=2,\ldots,m-1
\end{IEEEeqnarray*}

\begin{proposition}\label{p:Bt} The sets ${B}^t$ satisfy the following properties: \\
{\rm (d)} $\displaystyle B^t=\{\bar x\in F^n\,|\,p^t(\bar x)=\bar 0 \mbox{ and } \sum_{j,\ {\rm even}\, i}x^t_{i,j}=0 \}$;
\hfill{\rm\refstepcounter{equation}(\theequation)\label{eq:Bt}}\\
{\rm (d')} $|B^t|=|\Theta(A^t)|/2$;\\
{\rm (d'')} ${\rm Aut}(B^t)={\cal B}^t$.
\end{proposition}
\begin{proof}
(d) and (d') are straightforward.
For $t=1$ the claim (d'') trivially holds. Assume $t>1$.
Using Proposition~\ref{p:Theta}, we get
${\cal B}^t ={\rm Aut}(\{\bar r+\Theta(A^{t-1})\}_{\bar r\in V^t})$.
Moreover,
\begin{IEEEeqnarray*}{rCl}
{\rm Aut}(\{\bar r+\Theta(A^{t-1})\}_{\bar r\in V^t})
&=& {\rm Aut}\bigg(\bigcup_{\bar r\in V^t}(\bar r+\Theta(A^{t-1}))\bigg)
 \\ & = & {\rm Aut}(B^{t})
\end{IEEEeqnarray*}
because, as follows from Proposition \ref{p:At}(c), the sets $\bar r+\Theta(A^{t-1})$ are
connected components of
distance-$2$ graph of $B^{t}$.
\end{proof}

\begin{proposition}\label{p:AtBt}
Let $1\leq t \leq m-1$; then the following claims hold:\\
{\rm (a)} if $t<m-1$, then
${\cal A}^t=({\cal P}^t\rightthreetimes {\cal Q}^t)\rightthreetimes{\cal R}^t$
where
\begin{itemize}
 \item for groups ${\cal G}$ and ${\cal G}'$, the notation
  ${\cal G}\rightthreetimes{\cal G}'$ means a semidirect product with
  a normal subgroup ${\cal G}'$;
    \item
 ${\cal P}^t\simeq S_{2^{m-t}}$ is the subgroup of column permutations
 $\psi:(x^t_{i,j})_{i,j}\to (x^t_{i,\psi(j)})_{i,j}$;
    \item
 ${\cal Q}^t\simeq (S_{2^{t}})^{2^{m-t}}$
 is the set of collections of permutations in every column
 $(\phi_0,\ldots,\phi_{2^{m-t}-1}):(x^t_{i,j})_{i,j} \to (x^t_{\phi_j(i),j})_{i,j}$;
    \item
 ${\cal R}^t\simeq Z_2^{2^{m-t}(2^{t}-1)}$
 is the set of translations
 $\bar z+$, $\bar z\in \Theta(A^t)$;
\end{itemize}
{\rm (a')} ${\cal A}^{m-1}\simeq S_n\rightthreetimes Z_2^{n-1}$;\\
{\rm (b)} if $t<m-1$, then
${\cal B}^t=({\cal P}^t\rightthreetimes \widehat{\cal Q}^t)\rightthreetimes\widehat{\cal R}^t$
where
\begin{itemize}
    \item $\widehat{\cal Q}^t\simeq (S_2\rightthreetimes(S_{2^{t-1}})^2)^{2^{m-t}}$,
     $\widehat{\cal Q}^t \subset {\cal Q}^t$;
    \item
 $\widehat{\cal R}^t\simeq Z_2^{2^{m-t}(2^{t}-1)-1}$
 is the set of translations $\tau_{\bar z}$, $\bar z\in V^t+\Theta(A^{t-1})$ where
 $\tau_{\bar z}(\bar x)\eqdf \bar z+\bar x$.
\end{itemize}
{\rm (b')} ${\cal B}^{m-1}={\cal A}^{m-2}\rightthreetimes \{\tau_{\bar 0},\tau_{(11110...0)}\}$.
\end{proposition}
\begin{proof}
(a)
First we observe that ${\cal A}^t ={\rm Aut}(\Theta(A^t))$.
Since, by Proposition \ref{p:At}(c), $\Theta(A^t)$ is linear,
it holds ${\cal A}^t= {\cal O}^t \rightthreetimes {\cal R}^t$ where
${\cal O}^t\subset {\cal A}^t$ consists of coordinate permutations
and ${\cal R}^t\subset {\cal A}^t$ is a group of translations.

 It follows from Proposition \ref{p:At}(c) that ${\cal O}^t$ consists
of the permutations that do not break columns, i.\,e.,
an admissible permutation permutes columns and permutes elements in each column.

(a') follows from Proposition \ref{p:At}(c'').

(b) By Proposition~\ref{p:Bt}(d''), we have
${\cal B}^t = {\rm Aut}(B^{t})$.
Since $B^t$ is linear,
it holds ${\cal B}^t= \widehat{\cal O}^t \rightthreetimes \widehat{\cal R}^t$
where
$\widehat{\cal O}^t\subset {\cal B}^t$ is the coordinate permutation subgroup
and $\widehat{\cal R}^t\subset {\cal B}^t$ is the translation subgroup of ${\cal B}^t$.

Using Proposition~\ref{p:Bt}(d), we see that an arbitrary permutation from $\widehat{\cal O}^t$
does not break the columns of $(x^t_{i,j})_{i,j}$ and, moreover,
in each column the permutation does not change the parity of row-indexes
or changes the parity of all row-indexes.
(Indeed, in the case $t<m-1$ for any other coordinate permutation $\pi$ we can find a weight $2$ or weight $4$
word $\bar x$ such that $\bar x$ satisfies (\ref{eq:Bt}) but $\pi\bar x$ does not.)
It can be directly checked that all such permutations belong to ${\cal B}^t$.

(b') In the case $t=m-1$ the group ${\cal B}^t$ contains some other permutations
and this case can be easily calculated directly.
\end{proof}

\begin{corollary}\label{p:Dt}
$|{\cal D}^{m-1}|=n!/6((n/4)!)^4$.
 If $t<m-1$, then
$$
|{\cal D}^t|
=2\left( {\frac {2^t!}{2(2^{t-1}!)^2}} \right)^{2^{m-t}}
=2\left( {\frac 12 \big\lgroup {{2^t}\atop{2^{t-1}}} \big\rgroup} \right)^{2^{m-t}}
.$$
In particular,
$|{\cal D}^1|=2$,
$|{\cal D}^2|=2\cdot 3^{\frac{n}{4}}$,
$|{\cal D}^3|=2\cdot 35^{\frac{n}{8}}$,
$|{\cal D}^4|=2\cdot 6435^{\frac{n}{16}}$.
\end{corollary}

We say that an order-$t$ component $G$ is \emph{bold} {{if}} ${<}G{>}=B^t$ where
${<}G{>}$ means the affine span of $G$ (i.\,e., the minimal affine subspace
including $G$; if $G\ni \bar 0$, then the affine span coincides with the linear span).

The next proposition
helps us to see that all codes given by Construction \ref{con:3} are
pairwise different.
\begin{proposition}\label{p:bold}
Let $1\leq t\leq m-1$;
for each $\bar r\in V^t$ let $G_{\bar r}$ be a bold order-$(t-1)$ component and
$g_{\bar r}\in {\cal D}^{t-1}$. Put
$$ G \eqdf \bigcup_{\bar r\in V^t}(\bar r+g_{\bar r}(G_{\bar r})).$$
Then \\
{\rm (a)} $G$ is bold if and only if the collection
$\{g_{\bar r}\}_{\bar r\in V^t}$ is nondegenerate;\\
{\rm (b)} if $g'_{\bar r}\in {\cal D}^{t-1}$ and $G'_{\bar r}$ is an order-$(t-1)$  component for all $\bar r\in V^t$ (it is not necessary to assume that $G'_{\bar r}$ are bold), then
$$ G = \bigcup_{\bar r\in V^t}(\bar r+g'_{\bar r}(G'_{\bar r}))$$
implies $G'_{\bar r}=G_{\bar r}$ and $g'_{\bar r}=g_{\bar r}$ for all $\bar r\in V^t$.
\end{proposition}
\begin{proof}
(a) By the definition of bold component we have ${<}G_{\bar r}{>}=B^{t-1}$; thus
\begin{IEEEeqnarray*}{rCl}
{<}G{>}&=&{<}\bigcup_{\bar r\in V^t}(\bar r+g_{\bar r}(G_{\bar r})){>}
\\
&=&
{<}\bigcup_{\bar r\in V^t}(\bar r+g_{\bar r}({<}G_{\bar r}{>})){>}
\\
&=&
{<}\bigcup_{\bar r\in V^t}(\bar r+g_{\bar r}(B^{t-1})){>}.
\end{IEEEeqnarray*}
Since $g_{\bar r}(B^{t-1})$ is a half of $\Theta(A^{t-1})$,
the affine span ${<}G{>}$ coincides either with
$\bigcup_{\bar r\in V^t}(\bar r+\Theta(A^{t-1}))=B^t$ (i.\,e., $G$ is bold)
or with $\bigcup_{\bar r\in V^t}(\bar r+g_{\bar r}(B^{t-1}))$ ($G$ is not bold).
It is clear that the last case occurs if and only if
the sets $g_{\bar r}(B^{t-1})$, $\bar r\in V^t$, are translations of each other
(i.\,e., $g_{\bar r}$ have a common coordinate permutation)
and the translation vectors compose an affine function on $V^t$.


(b) It suffices to show that
for arbitrary $g,g'\in{\cal D}^{t-1}$ and bold components $G_0$, $G'_0$ of order $t-1$
the inequality $g'\neq g$ implies $g'(G'_0)\neq g(G_0)$.
This holds because, by the definitions of ${\cal D}^{t-1}$ and bold components and
the fact that ${\cal B}^{t-1}={\rm Aut}(B^{t-1})$ (Proposition~\ref{p:Bt}(d'')),
$g'\neq g$ implies $g'({<}G'_0{>})\neq g({<}G_0{>})$.
\end{proof}

\begin{proposition}\label{p:DtVt} If $1\leq t<m-1$, then
the number of degenerate collections $\{g_{\bar r}\in {\cal D}^t\}_{\bar r\in V^{t+1}}$
is $|{\cal D}^t|\cdot|V^{t+1}|$.
\end{proposition}
\begin{proof}
Assume $1\leq t<m-1$. As follows from Proposition~\ref{p:AtBt} and the fact that
$\Theta(A^t)/B^t=2$ (Proposition~\ref{p:Bt}(d')),
for each coordinate permutation $\pi$ there are $2$ or $0$ elements $\bar v$ such that
the automorphism $\bar v + \pi(\cdot)$ belongs to ${\cal D}^t$.
Thus we have:

1) The number of different coordinate permutations in ${\cal D}^t$ is $|{\cal D}^t|/2$.

2) For each admissible coordinate permutation $\pi$
the number of collections
$\{{\bar v_{\bar r}}+\pi(\cdot)\}_{\bar r\in V^{t+1}}$ of
automorphisms from ${\cal D}^t$ such that the set
$\{\bar r + \bar v_{\bar r} \,|\, {\bar r\in V^{t+1}} \}$
is an affine subspace equals
the number of two-value functions $f:V^{t+1}\to\{\bar v_1,\bar v_2\}$
satisfying
$f(\bar r_1)+f(\bar r_2)+f(\bar r_3)=f(\bar r_1+\bar r_2+\bar r_3)$
for any $\bar r_1,\bar r_2,\bar r_3\in V^{t+1}$,
i.\,e.,
 the number $2|V^{t+1}|$ of affine
$\{0,1\}$-value functions on $V^{t+1}$.

 By the definition of degenerate collection, the proposition is proved.
\end{proof}
\section{ A Lower Bound on the Number of $1$-Perfect Codes}
Denote by $\widetilde K_{\scriptscriptstyle LA}(n)$ the number of different
extended $1$-perfect codes given by Construction \ref{con:3}.

\begin{theorem}\label{th:Kla}
The extended $1$-perfect codes from Construction \ref{con:3}
are pairwise different. The number
of such codes equals
 \newcommand{\KKNK}{{\textstyle
    \big\lgroup\hspace{-0.7ex}  {{k} \atop {k{/}2}}  \hspace{-0.7ex}\big\rgroup
    }^{\hspace{-0.8ex}\frac{n}{k}}} 
 \newcommand{\DmMinusOne}{\frac{n!}{6\left({\frac n4}\mbox{\rm\large !}\right)^{4}}}

$$\begin{array}{l}
\widetilde K_{\scriptscriptstyle LA}(n)
=
|{\cal D}^{m-1}|\hspace{-0.5ex}\displaystyle\prod_{t=1}^{m-2}\hspace{-0.8ex}\left(
|{\cal D}^t|^{|V_{t+1}|}\hspace{-0.3ex}-\hspace{-0.3ex}|{\cal
D}^t|\hspace{-0.3ex}\cdot\hspace{-0.3ex}{|V_{t+1}|}\right)^{|V_{t+2}|\cdot\ldots \cdot|V_{m-1}|}
\nonumber\\
\qquad=
\DmMinusOne
\displaystyle\prod_{k=2,4,8,\ldots,\frac n4}
\Bigg(
\left( 2\cdot 2^{\vphantom{i}^{-\frac{n}{k}}} \KKNK \right)^
{\EExp[-0.5ex]{2k}}
\\
\hspace{30ex}-
\KKNK
\cdot{2^{\vphantom{i}^{-\frac{n}{2k}}}}
\Bigg)^{\EExpLog[-0.5ex]{2k}}
\end{array}$$
In particular, $\widetilde K_{\scriptscriptstyle LA}(16)=15 692 092 416 000 000$,
$\widetilde K_{\scriptscriptstyle LA}(32)\approx 2^{2363.79}$.
The following is the asymptotic formula for $\widetilde K_{\scriptscriptstyle LA}(n)$:
\begin{IEEEeqnarray}{rCl}
\widetilde K_{\scriptscriptstyle LA}(n)&\sim&|{\cal D}^{m-1}|\prod_{t=1}^{m-2}
|{\cal D}^t|^{|V_{t+1}|\cdot|V_{t+2}|\cdot\ldots \cdot|V_{m-1}|}
\nonumber\\&=&
\DmMinusOne
\prod_{k=2,4,8,\ldots,\frac n4}
\left(
{2}\cdot{2^{\vphantom{i}^{-\frac{n}{k}}}}
\KKNK
\right)^{\EExpLog[-0.5ex]{k}}
\label{eq:4}\\
&=&
                        \ExpHam{2}
\cdot
(   3^{\Eexp{4}}  \cdot \ExpHam{4} )
\nonumber\\&&
\qquad\qquad\,\,\,{}\cdot
(  35^{\Eexp{8}}  \cdot \ExpHam{8} )
\nonumber\\&&
\qquad\quad\,{}\cdot
(6435^{\Eexp{16}} \cdot \ExpHam{16} )
\nonumber\\&&
\,\cdot
\dots\cdot
\left(\left(
{\textstyle\frac 12{\big\lgroup\hspace{-0.7ex}}{{n/4}\atop{n/8}}
 {\hspace{-0.7ex}\big\rgroup}
}
\right)^{2^3}\cdot 2^{2^1}\right)
\cdot
\DmMinusOne
\nonumber
\end{IEEEeqnarray}
\end{theorem}
\begin{proof}
The number of ways to define an extended $1$-perfect code
using formulas (\ref{eq:2})
with restrictions of Construction~\ref{con:3}
can be easily calculated by Corollary~\ref{p:Dt} and
Proposition~\ref{p:DtVt}.
Proposition~\ref{p:bold} guarantees that
different local automorphisms give different codes.
\end{proof}
Since there is a one-to-one correspondence (deleting the last symbol) between extended
$1$-perfect codes and {\it $1$-perfect codes}, we have the following:

\begin{theorem}[A lower bound]
The number $B(n-1)$ of
$1$-perfect binary codes of length $n-1=2^m-1$
satisfies
\begin{equation}\label{eq:lbLA}
B(n-1)\geq \widetilde K_{\scriptscriptstyle LA}(n)
\end{equation}
where the exact expression and the asymptotic form for
$\widetilde K_{\scriptscriptstyle LA}(n)$
are given in
Theorem \ref{th:Kla}.
\end{theorem}

As we can see, the previous lower bound
$\ExpHam{2}\cdot(   3^{\Eexp{4}}  \cdot \ExpHam{4} )$
\cite{Kro:quas_l_b} consists of two multipliers ($k=2,4$)
of (\ref{eq:4}).


\begin{conjecture}\label{cj:} The lower bound {\rm (\ref{eq:lbLA})} is asymptotically tight,
i.\,e., (\ref{eq:4}) is the asymptotic number of $1$-perfect binary codes of length $n-1=2^m-1$.
\end{conjecture}

This conjecture is supported by our knowledge about $1$-perfect codes of small ranks, i.\,e., of rank $+1$ and of rank $+2$.
The \emph{rank} of the code is the dimension of its affine span;
we say that a $1$-perfect code of length $n-1$ is of rank $+p$
{{if}} its rank is $r_H+p$
where $r_H$ is the dimension of the linear $1$-perfect code (Hamming code) of corresponding length.
(The notion `affine span' means the same as `linear span'
if the code contains $\bar 0$, but the affine span is invariant for the code translations.)
We know that the LA construction gives almost all codes of rank $+1$
and almost all codes of rank $+2$ (and, of course, some other codes).
Moreover, if the affine span is fixed, then the number of $1$-perfect codes of rank $+1$ equals asymptotically
$$ \ExpHam{2} $$
and of rank $+2$,
$$
\ExpHam{2}\cdot
(   3^{\Eexp{4}}\cdot
    \ExpHam{4} )
    $$
(if we do not fix the affine span of code, then these values must be multiplied by the indexes
$n!/2^{n/2}(\frac{n}{2}-1)(\frac{n}{2}-2)\ldots(\frac{n}{2}-\frac{n}{4})$ and
$n!/24^{n/4}(\frac{n}{4}-1)(\frac{n}{4}-2)\ldots(\frac{n}{4}-\frac{n}{8})$ respectively).
This knowledge comes from the representation of $1$-perfect
binary codes of rank $+1$ and $+2$ \cite{AvgHedSol:class}
and the asymptotic number $3^{n+1}2^{2^n +1}(1+o(1))$
of $n$-ary quasigroups of order $4$ \cite{KroPot:Lyap,PotKro:asymp}.

\begin{remark}
All the codes given by Construction \ref{con:3} have the rank deficiency ${\rm RD}=2$
(the maximum rank of extended $1$-perfect binary codes of length $n\geq 16$ equals $n-1$, see \cite{EV:94}; so, the rank deficiency is defined as ${\rm RD}(C)\eqdf (n-1)-{\rm rank}(C)$)
and
(as follows from the bound
$\dim({\rm kernel}(C))\geq 2^{{\rm RD}(C)}$
for binary $1$-perfect codes of rank at least $+2$,
see \cite[Corollary 2.6]{PheVil:RankKernel})
the dimension of kernel at least $4$,
where ${\rm kernel}(C)\eqdf \{\bar k\,|\,C+\bar k=C\}$.
The last fact means that the construction gives at least
$$ \frac{\widetilde K_{\scriptscriptstyle LA}(n)}{n!2^{n-5}}$$
nonequivalent extended $1$-perfect binary codes of length $n$ and
$$ \frac{\widetilde K_{\scriptscriptstyle LA}(n)}{(n-1)!2^{n-5}}$$
nonequivalent $1$-perfect binary codes of length $n-1$,
where $n!2^{n-1}$ is the number of isometries of $\F$
and $(n-1)!2^{n-1}$ is the number of isometries of $F^{n-1}$.
\end{remark}

Yes, Conjecture~\ref{cj:} implies that almost all (extended) $1$-perfect codes have the
rank $n-3$, which is not full ($n-1$ for extended $1$-perfect codes of length $n$),
and even not fore-full ($n-2$).
This lacks support from the length-$16$ codes,
see \cite{ZZ:2006:l16r14} (the most part of codes has rank $14=n-2$, but the number of full-rank extended $1$-perfect codes of length $16$ is small indeed, see. \cite{OstPot:16},\cite{OPP:16}),
but the LA construction has not ``gathered power'' when $n=16$ ($m=4$).
Indeed, for $m=4$, among the three multipliers of (\ref{eq:4}), the first one ($t=1$)
is almost the same as the second ($t=2$), and the multiplier $n!/6((n/4)!)^4$
is the largest, while asymptotically the first multiplier is the most powerful one.
On the other hand, the fact that almost all codes have not full rank
can be expected. For example, this holds for $4$-ary distance $2$ MDS codes
($n$-ary qua\-si\-groups of order $4$, see \cite{PotKro:asymp,KroPot:n-ary});
the rank (over $Z_2^2$) has three different values for these codes
(rank $n-1$ for linear codes of length $n$, rank $n-\frac{1}{2}=\log_4|Z_2^{2n-1}|$ for `semilinear' codes, and rank $n$),
but the class with the middle rank value is the most powerful.
It is also notable that the number of nonequivalent order $16$
Steiner quadruple systems of full rank $15$ is smaller than of rank $14$ \cite{KOP}.
\pagebreak


\label{sec:biblio}
\providecommand\href[2]{#2} \providecommand\url[1]{\href{#1}{#1}}

\end{document}


\def\hrefDAOR#1{\href{http://math.nsc.ru/conference/DAOR'04/main.html}{#1}}
\renewcommand\thefootnote{}
\title{О числе $1$-со\-вер\-шен\-ных двоичных кодов: нижняя оценка%
\renewcommand\thefootnote{}
\thanks{
Это авторский перевод заметки в
\href{http://ieeexplore.ieee.org/xpl/RecentIssue.jsp?punumber=18}{IEEE Transactions on Information Theory}
54(4) 2008, 1760-1765,
DOI
\href{http://dx.doi.org/10.1109/TIT.2008.917692}{10.1109/TIT.2008.917692},
\copyright 2008 IEEE.}
\thanks{Материал настоящей заметки частично докладывался
на 10й Международной Конференции по Алгебраической и Комбинаторной Теории Кодирования
ACCT-10, Звенигород, Россия, сентябрь 2006, 
и, в краткой форме,
на Всероссийской конференции
``\hrefDAOR{Дискретный анализ} \hrefDAOR{и исследование операций}''
\hrefDAOR{{DAOR}'2004},
Новосибирск, Россия, июнь-июль 2004, \href{http://www.math.nsc.ru/conference/DAOR'04/daor4.pdf}{с.95}.}
}
\def\today{}
\author{\href{http://arxiv.org/find/grp_math/1/au:+Krotov_Denis/0/1/0/all/0/1}{Д.~С.~КРОТОВ},
\href{http://arxiv.org/find/grp_math/1/au:+Avgustinovich_Sergey/0/1/0/all/0/1}{С.~В.~АВГУСТИНОВИЧ}
\renewcommand\thefootnote{}
\thanks{Адрес авторов: Институт математики им. С. Л. Соболева СО РАН,
проспект Академика Коптюга 4, Новосибирск, 630090, Россия
(e-mail: krotov@math.nsc.ru, avgust@math.nsc.ru)}
 }
\maketitle

\begin{abstract}
Предложена конструкция $1$-со\-вер\-шен\-ных двоичных кодов,
дающая новую оценку снизу числа таких кодов.
Мы предполагаем, что эта оценка асимптотически точна.

\emph{Ключевые слова} --- {автоморфизм, оценка снизу, совершенный двоичный код}
\end{abstract}

\newcommand{\EExpLog}[2][0ex]{{{\raisebox{#1}{$2$}}^{\frac{n}{#2}-{\log}\frac{n}{#2}-1}}}
\newcommand{\EExp}[2][0ex]{{{\raisebox{#1}{$2$}}^{\frac{n}{#2}-1}}}
\newcommand{\EexpLog}[2][-0.2ex]{{{\raisebox{#1}{$\scriptstyle 2$}}^{\frac{n}{#2}-{\log}\frac{n}{#2}-1}}}
\newcommand{\Eexp}[2][-0.2ex]{{{\raisebox{#1}{$\scriptstyle 2$}}^{\frac{n}{#2}-1}}}
\newcommand{\ExpHam}[2][0ex]{2^{\EexpLog[-0.2ex]{#2}}}

\section{ Введение}
Работа посвящена проблеме перечисления $1$-со\-вер\-шен\-ных двоичных кодов.
Такие коды, как и любые оптимальные коды,
являются экстремальными объектами теории корректирующих кодов.
В то же время, совершенные коды --- это особый вид комбинаторных конфигураций.
Конструкция $1$-со\-вер\-шен\-ных двоичных кодов, предложенная в настоящей работе,
дает самый мощный известный класс таких кодов и приводит к оценке снизу их числа.

Первая известная конструкция \cite{LBru:Vas:nongroup_perfect}
нелинейных $1$-со\-вер\-шен\-ных двоичных кодов дает оценку
\begin{eqnarray} \nonumber
B(n-1)&\geq&\ExpHam{2}\cdot\ExpHam{4}\cdot\ExpHam{8}\cdot      \dots
\end{eqnarray}
числа $B(n-1)$ таких кодов длины $n-1=2^m-1$ (здесь и далее логарифм имеет основание $2$).
Эта оценка улучшалась в \cite{LBru:AvgSol:alpha-comp} и \cite{LBru:Mal:lb};
некоторые идеи, предложенные  в этих статьях, используются и в настоящей работе.
Наилучшая известная ранее оценка снизу \cite{LBru:Kro:quas_l_b} имеет вид
\begin{eqnarray} \label{LBru:eq:0}
B(n-1)&\geq&\ExpHam{2}\cdot 3^{\Eexp{4}}\cdot \ExpHam{4}.
\end{eqnarray}
Результат \cite{LBru:Kro:quas_l_b} сформулирован в терминах частного случая \cite{LBru:Phelps84}
обобщенной каскадной конструкции (см., напр., \cite{LBru:ZinLob:2000}),
который позволяет строить $1$-со\-вер\-шен\-ные двоичные коды из $q$-ичных МДР-кодов с расстоянием $2$
(для нижней оценки полезен лишь случай $q=4$),
или $n$-арных квазигрупп (порядка $q=4$).
Нижняя оценка числа $n$-арных квазигрупп порядка $4$, данная в \cite{LBru:Kro:quas_l_b},
асимптотически точна \cite{LBru:KroPot:Lyap},\cite{LBru:PotKro:asymp},
следовательно, этот способ оценки числа $1$-со\-вер\-шен\-ных двоичных кодов исчерпан.

Наилучшая известная оценка сверху \cite{LBru:Avg:ub} числа $1$-со\-вер\-шен\-ных двоичных кодов
имеет вид
$
2^{2^{n-(3/2){\log} n+\log{\log} (e n)}}.
$

Метод локальных автоморфизмов, предложенный в настоящей работе,
является дальнейшим развитием методов
 \cite{LBru:Vas:nongroup_perfect},\cite{LBru:AvgSol:alpha-comp},\cite{LBru:Mal:lb},\cite{LBru:Kro:quas_l_b}.
Поскольку имеется взаимнооднозначное соответствие между $1$-со\-вер\-шен\-ными двоичными кодами
и расширенными $1$-со\-вер\-шен\-ными двоичными кодами, результаты сформулированы в терминах последних.

\section{ Предварительные сведения}
Пусть $F^n$ ($\F$ / $\FF$) --- множество двоичных $n$-слов
(с четным / нечетным числом единиц, соответственно) с метрикой Хемминга $d(\cdot,\cdot)$
и покоординатным сложением по модулю~$2$.
Для слова $\bar x \in F^n$ определим его \emph{вес} $wt(\bar x) \eqdf d(\bar x,\bar 0)$.
 {\it Окрестностью} множества $S\subseteq F^n$ будем считать множество
$\Omega(S)\eqdf \bigcup_{\bar x\in S}\Omega(\bar x)$, где
$\Omega(\bar x)\eqdf \{\bar y\in F^n \,|\,  d(\bar y,\bar x)=1\}$.

Множество $C\subseteq F^n$ называется \emph{кодом с расстоянием $d$} (\emph{длины $n$}),
если расстояние Хемминга между любыми двумя различными словами из $C$ не меньше $d$.
 {\it Расширенный $1$-совершенный код\,} --- такое множество $C\subseteq \F$,
что окрестности слов из $C$ попарно не пересекаются и $\Omega(C)=\FF$.
Как следует из определения, $C$ суть код с расстоянием $4$ мощности
$|C|=|\FF|/n=2^{n-{\log} n-1}$ и $n$ --- степень двойки.
С другой стороны, любой код $C\subset \F$ с расстоянием $4$ мощности $|C|=2^{n-{\log} n-1}$,
очевидно, суть расширенный $1$-совершенный код.
Далее мы считаем, что $n=2^m\geq 16$.

Следующие формулы определяют некие вспомогательные множества $V^t,A^t\subset\F$
и дают представление расширенного кода Хемминга:
\begin{eqnarray}
V^t & \eqdf & \{(\bar{v},\bar v,0,\ldots ,0)\in F^n\,|\,\bar v\in \F[2^{m-t}]\}, \label{LBru:eq:Vt}
\\ \nonumber&&
 t=1,\ldots,m-1
\\
 A^1 & \eqdf & V^1
\nonumber \\
A^{t} & \eqdf & V^t + A^{t-1} = \bigcup_{\bar r\in V^{t}}(\bar r+A^{t-1}), \label{LBru:eq:At}
\\ \nonumber&&
 t=2,\ldots,m-1
 \\
H & \eqdf & A^{m-1}. \label{LBru:eq:1}
\end{eqnarray}
Непосредственно проверяется следующий факт:
\begin{proposition}
Множество $H$, определяемое формулами {\rm (\ref{LBru:eq:Vt})-(\ref{LBru:eq:1})},
является линейным расширенным $1$-совершенным кодом,
т.\,е. расширенным кодом Хемминга (см., напр., {\rm\cite[\S1.7]{LBru:MWS}}),
который, как известно, есть единственный линейный расширенный $1$-совершенный код,
с точностью до перестановки координат.
\end{proposition}

Будем говорить, что множество $G\subset \F$ является \emph{компонентой порядка $t\in\{1,\ldots,m-1\}$},
если $|G|=|A^t|$ и $\Omega(G)=\Omega(A^{t})$.
Пусть $t\in\{1,\ldots,m-1\}$ и $\bar \mu\in F^{2^{m-t}}\times {0}^{n-2^{m-t}}$;
будем говорить, что множество $M\subset \F$
является \emph{$\bar \mu$-компонентой} (порядка $t$),
если $M=G+\bar \mu$ для некоторой компоненты $G$ порядка $t$.

Пусть ${\rm Aut}(F^n)$ обозначает группу изометрий $F^n$
(которая совпадает с группой автоморфизмов графа расстояний $1$ на вершинах $F^n$).
Известно, что каждая изометрия $g\in {\rm Aut}(F^n)$ имеет единственное представление
 $g(\cdot) = \bar v+\pi(\cdot)$, где $\bar v$ --- \emph{вектор сдвига}
 из $F^n$ и $\pi$ перестановка координат.
Если $\pi=Id$, т.\,е. $g(\cdot) = \bar v+\cdot$, то изометрия $g$ называется \emph{сдвигом}.
Для множества $S\subseteq F^n$ обозначим через ${\rm Aut}({S})$
группу изометрий $g\in{\rm Aut}(F^n)$ таких, что $g(S)=S$.
Для набора ${\mathbf{S}}=\{S_1,\ldots ,S_l\}$ подмножеств $F^n$
обозначим через ${\rm Aut}(\mathbf{S})$ группу изометрий $g\in{\rm Aut}(F^n)$ таких, что
для каждого $S\in\mathbf{S}$ множество $g(S)$ также принадлежит $\mathbf{S}$.
Далее мы будем использовать каллиграфические буквы для обозначения подгрупп ${\rm Aut}(F^n)$.
Положим $${\cal A}^t \eqdf {\rm Aut}(\Omega(A^{t})),$$
где $A^t$ определяется формулами (\ref{LBru:eq:Vt})-(\ref{LBru:eq:At}).

\begin{proposition}\label{LBru:pro:1} {\rm (a)} Пусть $t\in\{2,\ldots,m-1\}$ и
пусть для каждого $\bar\mu\in V^t$ множество
$B_{\bar\mu}$ есть ${\bar\mu}$-компонента порядка $t-1$.
Тогда множество $B=\bigcup_{\bar\mu\in V^t} B_{\bar\mu}$ является компонентой порядка $t$.

{\rm (b)} Если $B$ --- компонента порядка $t$ и $g^t\in {\cal A}^t$,
то $g^t(B)$ --- также компонента порядка $t$.

{\rm (c)} Компонента порядка $m-1$ является расширенным $1$-совершенным кодом.
\end{proposition}
\begin{proof}
П.(a) следует из (\ref{LBru:eq:At}), определения компоненты порядка $t$
и равенства $\Omega(\bar\mu+A^{t-1})=\Omega(B_{\bar\mu})$,
которое выполняется по определению $\bar\mu$-компоненты порядка $t-1$.

Пп.(b) и (c) следуют прямо из определений.
\end{proof}

\section{ ЛА конструкция расширенных $1$-совершенных кодов}
Предложения~\ref{LBru:pro:1} достаточно, чтобы удостовериться, что конструкция,
приведенная ниже, строит расширенные $1$-совершенные коды.
Идея конструкции --- стартуя с кода Хемминга, применять изометрии $F^n$
к частям кода таким образом, чтобы окрестности этих частей не менялись.
Мы называем такие изометрии \emph{локальными автоморфизмами}; локальный автоморфизм
действует на часть кода, не меняя ее окрестности.
На первом этапе мы берем в качестве таких частей компоненты порядка $1$,
на втором --- компоненты порядка $2$, и так далее.
На последнем этапе мы ``крутим'' весь код.

\begin{constr}[ЛА --- локальные автоморфизмы] \label{LBru:con:1}
Пусть для каждого целого $t \in \{2,\ldots, m\}$ и для каждого слова
${\bar r_{i}\in V^{i}}$, $i=t,\ldots,m-1$
определен локальный автоморфизм $g_{\bar r_t,\ldots ,\bar r_{m-1}}\in {\cal A}^{t-1}$;
в частности, $g\in {\cal A}^{m-1}$.
Тогда (как следует индукцией по $t$ из Предложения~\ref{LBru:pro:1})
множество $C$, представленное следующими формулами, является расширенным $1$-совершенным кодом.
\begin{eqnarray} 
A^1_{\bar r_{2},\ldots ,\bar r_{m-1}} & \eqdf & V^1,\qquad \bar r_{i}\in V^i
\nonumber\\
A^{t}_{\bar r_{t+1},\ldots ,\bar r_{m-1}} & \eqdf & \bigcup_{\bar r_{t}\in V^{t}}
\left(\bar r_t+g_{\bar r_t,\ldots ,\bar r_{m-1}}(A^{t-1}_{\bar r_t,\ldots ,\bar r_{m-1}})\right),
\nonumber \\ \nonumber&&
 \qquad\qquad t=2,\ldots,m-1
\label{LBru:eq:2}
\\
C & \eqdf & g(A^{m-1}).
\end{eqnarray}
\end{constr}

В Конструкции~\ref{LBru:con:1} каждый код получается более чем одним способом.
Чтобы оценить число кодов, получаемых таким образом, нужны более строгие ограничения на локальные автоморфизмы.

Пусть
\begin{eqnarray*}
{\cal B}^1 &\eqdf& {\rm Aut}(A^{1})
\\
{\cal B}^t &\eqdf& {\rm Aut}\big(\{\bar r+\Omega(A^{t-1})\}_{\bar r\in V^{t}}\big),\qquad t=2,\ldots ,m-1.
\end{eqnarray*}
Для каждого $t=1,\ldots ,m-1$ зафиксируем множество ${\cal D}^{t}$ представителей
смежных классов из ${\cal A}^{t}/{\cal B}^{t}$.
Более того, мы выберем представителей таким образом, чтобы выполнялось следующее условие:
для двух смежных классов $D_1,D_2\in{\cal A}^{t}/{\cal B}^{t}$ и их представителей
$d_1,d_2\in {\cal D}^{t}$, $d_1 \in D_1$, $d_2\in D_2$, равенство
$D_1=\tau D_2$ с некоторым сдвигом $\tau$ влечет $d_1=\tau d_2$
(это условие существенно в определении вырожденного набора
и в Предложении~\ref{LBru:p:bold} ниже).

По индукции можно показать, что

\begin{proposition}
Ограничения
$
g_{\bar r_t,\ldots ,\bar r_{m-1}}\in {\cal D}^{t-1}
$
не уменьшают множество кодов, представленных формулами {\rm (\ref{LBru:eq:2})}.
\end{proposition}
\begin{proof}
Пусть
$$G^t \eqdf \bigcup_{\bar r \in V^t}\left(\bar r+g_{\bar r}(G^{t-1}_{\bar r})\right),$$
где $G^{t-1}_{\bar r}$ --- компонента порядка $t-1$
и $g_{\bar r}\in{\cal A}^{t-1}$
для всех $\bar r\in V^t$.
Пусть $g\in {\cal A}^t$ и $g=dh$, где $d\in{\cal D}^t$
и $h\in {\cal B}^t$.

Мы утверждаем, что
\begin{equation}\label{LBru:eq:cl}
    g (G^t) = d (G'^t), \quad \mbox{где} \quad G'^t \eqdf \bigcup_{\bar q \in V^t}\left(\bar q+g'_{\bar q}(G^{t-1}_{\rho\bar q})\right)
\end{equation}
для некоторых $g'_{\bar q}\in{\cal A}^{t-1}$ и перестановки $\rho:V^t\to V^t$.
Действительно, по определению ${\cal B}^t$ для каждого $\bar r\in V^t$ имеем
$$
h(\bar r+\Omega(A^{t-1}))=\rho^{-1}\bar r+\Omega(A^{t-1}),
$$
где $\rho$ --- некоторая перестановка на $V^t$.
Таким образом, мы видим, что $h_{\bar r}(\cdot)\eqdf \rho^{-1}\bar r+h(\bar r+\cdot)$
принадлежит ${\cal A}^{t-1}$.
Далее, подставляя $\bar q\eqdf \rho^{-1}\bar r$ вместо $\bar r$,
мы видим, что (\ref{LBru:eq:cl}) выполняется с
$g'_{\bar q} \eqdf h_{\rho\bar q} g_{\rho\bar q}$.

Таким образом, используя (\ref{LBru:eq:cl}), мы можем шаг за шагом заменить операторы
$g_{...}\in {\cal A}^t$ на $d_{...}\in {\cal D}^t$,
стартуя с $t=m-1$ и заканчивая $t=1$.
\end{proof}

Следовательно, следующая конструкция дает то же множество кодов,
 что и Конструкция~\ref{LBru:con:1}.
\begin{constr}[ЛА, оценка сверху] \label{LBru:con:2}
Пусть для каждого целого $t \in \{2,\ldots, m\}$ и для любых слов
${\bar r_{i}\in V^{i}}$, $i=t,\ldots,m-1$,
мы имеем $g_{\bar r_t,\ldots ,\bar r_{m-1}}\in {\cal D}^{t-1}$;
в частности, $g\in {\cal D}^{m-1}$.
Тогда множество $C$, определяемое формулами {\rm (\ref{LBru:eq:2})}, является расширенным $1$-совершенным кодом.
\end{constr}

Как мы увидим ниже (Теорема~\ref{LBru:th:Kla}),
почти все ($n\to\infty$) коды, представленные Конструкцией~\ref{LBru:con:2},
имеют единственное представление, что дает хорошую оценку сверху
\begin{equation}\label{LBru:up}K_{\scriptscriptstyle LA}(n)\leq
|{\cal D}^{m-1}|\prod_{t=1}^{m-2}
|{\cal D}^t|^{|V_{t+1}|\cdot|V_{t+2}|\cdot\ldots \cdot|V_{m-1}|}
\end{equation}
числа $K_{\scriptscriptstyle LA}(n)$ различных расширенных
$1$-совершенных кодов длины $n$,
полученных методом локальных автоморфизмов (\emph{ЛА}), т.\,е.
Конструкцией~\ref{LBru:con:1} или~\ref{LBru:con:2}.
Чтобы показать, что число различных ЛА кодов близко к этому значению,
мы введем еще некоторые ограничения на $g_{\bar r_t,\ldots ,\bar r_{m-1}}$.

Пусть $L$ --- линейное подпространство $F^n$ и для каждого $\bar r\in L$ имеем
$g_{\bar r}\in {\cal D}^{t}$ и
$g_{\bar r}(\cdot)= \bar v_{\bar r}+\pi_{\bar r}(\cdot)$.
Будем говорить, что набор $\{g_{\bar r}\}_{\bar r\in L}$  \emph{вырожденный},
если выполнены следующие условия:

$\bullet$ перестановка $\pi_{\bar r}$ не зависит от $\bar r$,
т.\,е. $\pi_{\bar r}=\pi$ для всех $\bar r\in L$;

$\bullet$ множество $\{\bar r + \bar v_{\bar r} \,|\, {\bar r\in L} \}$
есть аффинное подпространство $F^n$.

В противном случае будем говорить, что набор $\{g_{\bar r}\}_{\bar r\in L}$ \emph{невырожденный}.

\begin{constr}[ЛА, оценка снизу] \label{LBru:con:3}
Дополнительно к условиям Конструкции~\ref{LBru:con:2} мы требуем, чтобы
набор $\mathbf{g}_{\bar r_{t+1},\ldots ,\bar r_{m-1}}=\{g_{\bar r_t,\ldots ,\bar r_{m-1}}\in {\cal D}^{t-1}\}_{\bar r_t\in V^t}$
был невырожденным для каждых
$t\in \{2,\ldots,{m-1}\}$, $\bar r_{t+1}\in  V^{t+1}$, \ldots, $\bar r_{m-1}\in  V^{m-1}$.
\end{constr}

\section{ Вычисления}

В этом разделе мы установим некторые факты касательно структуры компонент порядка $t$
и связанных с ними объектов, на которых базируется доказательство основного результата.
Для подмножества $G\subseteq \F$ положим
$$\Theta(G) \eqdf \{\bar x\in \F\,|\, \Omega(\bar x)\subseteq \Omega(G)\};$$
очевидно, $G\subseteq\Theta(G)$ и $\Omega(\Theta(G))=\Omega(G)$.
Следующий факт также следует прямо из определений:
\begin{proposition} \label{LBru:p:Theta} Для любых подмножеств $G,G'\subseteq \F$ равенства
$\Omega(G)=\Omega(G')$ и $\Theta(G)=\Theta(G')$ эквивалентны.
\end{proposition}
Для каждых $t=1,\ldots,m$ и $\bar x=(\bar x_0,\ldots,\bar x_{2^t-1})\in (F^{2^{m-t}})^{2^t}=F^n$
определим обобщенную проверку на четность
$$p^t(\bar x) \eqdf \sum_{i=0}^{2^t-1}\bar x_i.$$

\begin{proposition} \label{LBru:p:At}
Пусть $1\leq t \leq m-1$. Тогда выполнены следующие утверждения:\\
{\rm (a)} $p^t(\bar x)=\bar 0$ для всех $\bar x\in A^t;$\\
{\rm (a')} $|A^t|=2^{2^{m-t}(2^t-1)-t};$\\
{\rm (b)} $\Omega(A^t)=\{\bar x\in F^n\,|\,wt(p^t(\bar x))=1\};$\\
{\rm (b')} $|\Omega(A^t)|=2^{2^{m-t}(2^t-1)+m-t};$\\
{\rm (c)} если $t<m-1$, то $\Theta(A^t)=\{\bar x\in F^n\,|\,p^t(\bar x)=\bar 0\};$\\
{\rm (c')} если $t<m-1$, то $|\Theta(A^t)|=2^{2^{m-t}(2^t-1)};$\\
{\rm (c'')} $\Theta(A^{m-1})=\F$.
\end{proposition}
\begin{proof}
(a) и (a') следуют прямо из определения $A^t$.

(b') Поскольку $A^t$ имеет кодовое расстояние $4$, имеем $|\Omega(A^t)|=n|A^t|$.

(b) Из (a) следует, что $wt(p^t(\bar x))=1$ для всех $\bar x\in\Omega(A^t)$.
С другой стороны, согласно (b'), имеем $|\Omega(A^t)|=|\{\bar x\in F^n\,|\,wt(p^t(\bar x))=1\}|$.

(c) Согласно (b), из равенства $p^t(\bar x)=\bar 0$ следует $\bar x\in \Theta(A^t)$.
Допустим, что $p^t(\bar x)\neq\bar 0$. Если $t<m-1$, то найдется
$\bar y\in \Omega(\bar x)$ такое, что $wt(p^t(\bar y))>1$. Следовательно, $\bar x\not\in\Theta(A^t)$.

(c') следет из (c).

(c'') Из (b) и (b') следует, что $\Omega(A^{m-1})=\FF$, откуда $\Theta(A^{m-1})=\F$.
\end{proof}

Далее мы будем пользоваться представлением элементов из $F^n$ в виде массивов:
$$\bar x=(x^t_{0,0},...,x^t_{0,2^{m-t}-1},x^t_{1,0},
\ldots,
x^t_{2^{t}-1,2^{m-t}-1})=
(x^t_{i,j})_{i,j},$$
где индексы $i,j$ изменяются в лексикографическом порядке.
Т.\,е. для каждого $t=1,\ldots,m-1$ элемент $\bar x$ из $F^n$
можно представить как матрицу размера $2^{t}\times 2^{m-t}$
$$
\left(
\begin{array}{l@{\ }l@{\ }l@{\ \ }l}
x^t_{0,0} & x^t_{0,1} & \ldots & x^t_{0,2^{m-t}-1}\\
\ldots &\ldots &\ldots &\ldots \\
x^t_{2^{t}-1,0}& x^t_{2^{t}-1,1}& \ldots & x^t_{2^{t}-1,2^{m-t}-1}
\end{array}
\right)
$$
В этих терминах $p^t(\bar x)$ есть сумма строк матрицы $(x^t_{i,j})_{i,j}$.
Для дальнейших вычислений, определим множества
\begin{eqnarray*}
  {B}^1 & \eqdf & V^1  \\
  {B}^t & \eqdf & V^t + \Theta(A^{t-1}) 
  , \qquad t=2,\ldots,m-1
\end{eqnarray*}

\begin{proposition}\label{LBru:p:Bt} Множества ${B}^t$ обладают следующими свойствами: \\
{\rm (d)} $\displaystyle B^t=\{\bar x\in F^n\,|\,p^t(\bar x)=\bar 0 \mbox{ и }
           \sum_{j,\ \mbox{\rm\footnotesize чётн.}\, i}x^t_{i,j}=0 \}$;
\hfill{\rm\refstepcounter{equation}(\theequation)\label{LBru:eq:Bt}}\\
{\rm (d')} $|B^t|=|\Theta(A^t)|/2$;\\
{\rm (d'')} ${\rm Aut}(B^t)={\cal B}^t$.
\end{proposition}
\begin{proof}
Пп. (d) и (d') очевидны.
Для $t=1$ утверждение (d'') тривиално выполняется. Пусть $t>1$.
Пользуясь Предложением~\ref{LBru:p:Theta}, получаем
${\cal B}^t ={\rm Aut}(\{\bar r+\Theta(A^{t-1})\}_{\bar r\in V^t})$.
Более того,
\begin{eqnarray*}
{\rm Aut}(\{\bar r+\Theta(A^{t-1})\}_{\bar r\in V^t})
&=& {\rm Aut}\bigg(\bigcup_{\bar r\in V^t}(\bar r+\Theta(A^{t-1}))\bigg)
 \\ & = & {\rm Aut}(B^{t}),
\end{eqnarray*}
поскольку, как следует из Предложения~\ref{LBru:p:At}(c), множества $\bar r+\Theta(A^{t-1})$
являются компонентами связности графа расстояний $2$ множества $B^{t}$.
\end{proof}

\begin{proposition}\label{LBru:p:AtBt}
Пусть $1\leq t \leq m-1$, тогда верны следующие факты:\\
{\rm (a)} если $t<m-1$, то
${\cal A}^t=({\cal P}^t\rightthreetimes {\cal Q}^t)\rightthreetimes{\cal R}^t$,
где
\begin{itemize}
 \item для групп ${\cal G}$ and ${\cal G}'$ обозначение
  ${\cal G}\rightthreetimes{\cal G}'$ используется для их полупрямого произведения, где
  ${\cal G}'$ --- нормальная подгруппа;
    \item
 ${\cal P}^t\simeq S_{2^{m-t}}$ --- подгруппа перестановок столбцов
 $\psi:(x^t_{i,j})_{i,j}\to (x^t_{i,\psi(j)})_{i,j}$;
    \item
 ${\cal Q}^t\simeq (S_{2^{t}})^{2^{m-t}}$
 --- множество наборов перестановок в каждом столбце
 $(\phi_0,\ldots,\phi_{2^{m-t}-1}):(x^t_{i,j})_{i,j} \to (x^t_{\phi_j(i),j})_{i,j}$;
    \item
 ${\cal R}^t\simeq Z_2^{2^{m-t}(2^{t}-1)}$
 --- множество сдвигов
 $\bar z+$, $\bar z\in \Theta(A^t)$;
\end{itemize}
{\rm (a')} ${\cal A}^{m-1}\simeq S_n\rightthreetimes Z_2^{n-1}$;\\
{\rm (b)} если $t<m-1$, то
${\cal B}^t=({\cal P}^t\rightthreetimes \widehat{\cal Q}^t)\rightthreetimes\widehat{\cal R}^t,$
где
\begin{itemize}
    \item $\widehat{\cal Q}^t\simeq (S_2\rightthreetimes(S_{2^{t-1}})^2)^{2^{m-t}}$,
     $\widehat{\cal Q}^t \subset {\cal Q}^t$;
    \item
 $\widehat{\cal R}^t\simeq Z_2^{2^{m-t}(2^{t}-1)-1}$
 --- множество сдвигов $\tau_{\bar z}$, $\bar z\in V^t+\Theta(A^{t-1}),$ где
 $\tau_{\bar z}(\bar x)\eqdf \bar z+\bar x$.
\end{itemize}
{\rm (b')} ${\cal B}^{m-1}={\cal A}^{m-2}\rightthreetimes \{\tau_{\bar 0},\tau_{(11110...0)}\}$.
\end{proposition}
\begin{proof}
(a)
Заметим сначала, что ${\cal A}^t ={\rm Aut}(\Theta(A^t))$.
Поскольку по Предложению~\ref{LBru:p:At}(c) множество $\Theta(A^t)$ линейное,
выполнено равенство ${\cal A}^t= {\cal O}^t \rightthreetimes {\cal R}^t$, где
${\cal O}^t\subset {\cal A}^t$ состоит из перестановок координат
и ${\cal R}^t\subset {\cal A}^t$ --- группа сдвигов.

 Из Предложения~\ref{LBru:p:At}(c) следует, что ${\cal O}^t$ состоит
из перестановок, не меняющих разбиение на столбцы, т.\,е.
допустимая перестановка переставляет столбцы между собой и переставляет элементы внутри каждого столбца.

(a') следует из Предложения~\ref{LBru:p:At}(c'').

(b) По Предложению~\ref{LBru:p:Bt}(d'') имеем
${\cal B}^t = {\rm Aut}(B^{t})$.
Поскольку множество $B^t$ линейно,
выполнено ${\cal B}^t= \widehat{\cal O}^t \rightthreetimes \widehat{\cal R}^t$,
где
$\widehat{\cal O}^t\subset {\cal B}^t$ --- подгруппа перестановок координат
и $\widehat{\cal R}^t\subset {\cal B}^t$ --- подгруппа сдвигов группы ${\cal B}^t$.

Пользуясь Предложением~\ref{LBru:p:Bt}(d), мы видим, что произвольная перестановка из $\widehat{\cal O}^t$
не меняет разбиение на столбцы массива $(x^t_{i,j})_{i,j}$ и, более того,
в каждом столбце перестановка не меняет четность индексов строк
либо меняет четность одновременно всех индексов строк.
(Действительно, в случаe $t<m-1$ для любой другой перестановки $\pi$
мы можем найти слово $\bar x$ веса $2$ или $4$ такое, что
$\bar x$ удовлетворяет (\ref{LBru:eq:Bt}), а $\pi\bar x$ --- нет.)
Легко убедиться, что все такие перестановки принадлежат ${\cal B}^t$.

(b') В случае $t=m-1$ группа ${\cal B}^t$ содержит некоторые добавочные перестановки,
и этот случай легко проверяется непосредственно.
\end{proof}

\begin{corollary}\label{LBru:p:Dt}
$|{\cal D}^{m-1}|=n!/6((n/4)!)^4$.
 Если $t<m-1$, то
$$
|{\cal D}^t|
=2\left( {\frac {2^t!}{2(2^{t-1}!)^2}} \right)^{2^{m-t}}
=2\left( {\frac 12 \big\lgroup {{2^t}\atop{2^{t-1}}} \big\rgroup} \right)^{2^{m-t}}
.$$
В частности,
$|{\cal D}^1|=2$,
$|{\cal D}^2|=2\cdot 3^{\frac{n}{4}}$,
$|{\cal D}^3|=2\cdot 35^{\frac{n}{8}}$,
$|{\cal D}^4|=2\cdot 6435^{\frac{n}{16}}$.
\end{corollary}

Будем говорить, что компонента  $G$ порядка $t$ \emph{плотная}, {{если}} ${<}G{>}=B^t$, где
${<}G{>}$ означает аффинную оболочку множества $G$ (т.\,е. минимальное аффинное подпространство,
включающее $G$; в случае $G\ni \bar 0$ аффинная оболочка совпадает с линейной оболочкой).

Следующее предложение поясняет, почему все коды
из Конструкции~\ref{LBru:con:3} попарно различны.
\begin{proposition}\label{LBru:p:bold}
Пусть $1\leq t\leq m-1$.
Для каждого $\bar r\in V^t$ пусть $G_{\bar r}$ есть плотная компонента порядка $(t-1)$ и
$g_{\bar r}\in {\cal D}^{t-1}$. Положим
$$ G \eqdf \bigcup_{\bar r\in V^t}(\bar r+g_{\bar r}(G_{\bar r})).$$
Тогда \\
{\rm (a)} Компонента $G$ плотная если и только если набор
$\{g_{\bar r}\}_{\bar r\in V^t}$ невырожденный;\\
{\rm (b)} если $g'_{\bar r}\in {\cal D}^{t-1}$ и для каждого $\bar r\in V^t$ множество $G'_{\bar r}$ ---  компонента порядка $(t-1)$ (здесь мы не подразумеваем, что компонента $G'_{\bar r}$ обязательно плотная), то из равенства
$$ G = \bigcup_{\bar r\in V^t}(\bar r+g'_{\bar r}(G'_{\bar r}))$$
следует $G'_{\bar r}=G_{\bar r}$ и $g'_{\bar r}=g_{\bar r}$ для всех $\bar r\in V^t$.
\end{proposition}
\begin{proof}
(a) По определению плотной компоненты имеем ${<}G_{\bar r}{>}=B^{t-1}$, откуда
\begin{eqnarray*}
{<}G{>}&=&{<}\bigcup_{\bar r\in V^t}(\bar r+g_{\bar r}(G_{\bar r})){>}
\\
&=&
{<}\bigcup_{\bar r\in V^t}(\bar r+g_{\bar r}({<}G_{\bar r}{>})){>}
\\
&=&
{<}\bigcup_{\bar r\in V^t}(\bar r+g_{\bar r}(B^{t-1})){>}.
\end{eqnarray*}
Поскольку $g_{\bar r}(B^{t-1})$ --- половина от $\Theta(A^{t-1})$,
аффинная оболочка ${<}G{>}$ совпадает либо с
$\bigcup_{\bar r\in V^t}(\bar r+\Theta(A^{t-1}))=B^t$ (т.\,е. $G$ плотная),
либо с $\bigcup_{\bar r\in V^t}(\bar r+g_{\bar r}(B^{t-1}))$ ($G$ не плотная).
Ясно, что последний случай имеет место если и только если
множества $g_{\bar r}(B^{t-1})$, $\bar r\in V^t$, --- сдвиги друг друга
(т.\,е. $g_{\bar r}$ имеют одну и ту же перестановку координат)
и векторы сдвигов образуют аффинную функцию на $V^t$.


(b) Достаточно показать, что для произвольных $g,g'\in{\cal D}^{t-1}$
и полных компонент $G_0$, $G'_0$ порядка $t-1$
из $g'\neq g$ следует $g'(G'_0)\neq g(G_0)$.
Это верно, поскольку из определений ${\cal D}^{t-1}$ и полных компонент и равенства
${\cal B}^{t-1}={\rm Aut}(B^{t-1})$ (Предложение~\ref{LBru:p:Bt}(d'')) имеем:
 $g'\neq g$ подразумевает $g'({<}G'_0{>})\neq g({<}G_0{>})$.
\end{proof}

\begin{proposition}\label{LBru:p:DtVt} Если $1\leq t<m-1$, то
число вырожденных наборов $\{g_{\bar r}\in {\cal D}^t\}_{\bar r\in V^{t+1}}$
равно $|{\cal D}^t|\cdot|V^{t+1}|$.
\end{proposition}
\begin{proof}
Пусть $1\leq t<m-1$. Как следует из Предложения~\ref{LBru:p:AtBt} и равенства $\Theta(A^t)/B^t=2$ (Предложение~\ref{LBru:p:Bt}(d')),
для каждой перестановки координат $\pi$ имеется $2$ или $0$ элементов $\bar v$ таких, что
автоморфизм $\bar v + \pi(\cdot)$ принадлежит ${\cal D}^t$.
Отсюда имеем следующее:

1) Число различных перестановок координат в ${\cal D}^t$ равно $|{\cal D}^t|/2$.

2) Для каждой допустимой перестановки координат $\pi$
число таких наборов
$\{{\bar v_{\bar r}}+\pi(\cdot)\}_{\bar r\in V^{t+1}}$  автоморфизмов из ${\cal D}^t$,
что множество
$\{\bar r + \bar v_{\bar r} \,|\, {\bar r\in V^{t+1}} \}$
есть аффинное подпространство, равно числу двузначных функций $f:V^{t+1}\to\{\bar v_1,\bar v_2\}$,
удовлетворяющих соотношению
$f(\bar r_1)+f(\bar r_2)+f(\bar r_3)=f(\bar r_1+\bar r_2+\bar r_3)$
для любых $\bar r_1,\bar r_2,\bar r_3\in V^{t+1}$,
т.\,е.
 числу $2|V^{t+1}|$ аффинных
$\{0,1\}$-значных функций на $V^{t+1}$.

По определению вырожденного набора, предложение доказано.
\end{proof}
\section{ Нижняя оценка \\\mbox{}\hfill числа $1$\mbox{-}совершенных кодов}
Обозначим через $\widetilde K_{\scriptscriptstyle LA}(n)$ число различных
расширенных $1$-совершенных кодов, представляемых Конструкцией~\ref{LBru:con:3}.

\begin{theorem}\label{LBru:th:Kla}
Расширенные $1$-совершенные коды из Конструкции~\ref{LBru:con:3}
попарно различны. Число таких кодов равно
 \newcommand{\KKNK}{{\textstyle
    \big\lgroup\hspace{-0.7ex}  {{k} \atop {k{/}2}}  \hspace{-0.7ex}\big\rgroup
    }^{\hspace{-0.8ex}\frac{n}{k}}} 
 \newcommand{\DmMinusOne}{\frac{n!}{6\left({\frac n4}\mbox{\rm\large !}\right)^{4}}}

$$\begin{array}{l}
\widetilde K_{\scriptscriptstyle LA}(n)
=
|{\cal D}^{m-1}|\hspace{-0.5ex}\displaystyle\prod_{t=1}^{m-2}\hspace{-0.8ex}\left(
|{\cal D}^t|^{|V_{t+1}|}\hspace{-0.3ex}-\hspace{-0.3ex}|{\cal
D}^t|\hspace{-0.3ex}\cdot\hspace{-0.3ex}{|V_{t+1}|}\right)^{|V_{t+2}|\cdot\ldots \cdot|V_{m-1}|}
\nonumber\\
\qquad=
\DmMinusOne
\displaystyle\prod_{k=2,4,8,\ldots,\frac n4}
\Bigg(
\left( 2\cdot 2^{\vphantom{i}^{-\frac{n}{k}}} \KKNK \right)^
{\EExp[-0.5ex]{2k}}
\\
\hspace{30ex}-
\KKNK
\cdot{2^{\vphantom{i}^{-\frac{n}{2k}}}}
\Bigg)^{\EExpLog[-0.5ex]{2k}}
\end{array}$$
В частности, $\widetilde K_{\scriptscriptstyle LA}(16)=15 692 092 416 000 000$,
$\widetilde K_{\scriptscriptstyle LA}(32)\approx 2^{2363.79}$.
Следующая формула дает асимптотику числа $\widetilde K_{\scriptscriptstyle LA}(n)$:
\begin{eqnarray}
\widetilde K_{\scriptscriptstyle LA}(n)&{\sim}&|{\cal D}^{m-1}|\prod_{t=1}^{m-2}
|{\cal D}^t|^{|V_{t+1}|\cdot|V_{t+2}|\cdot\ldots \cdot|V_{m-1}|}
\nonumber\\&{=}&
\DmMinusOne
\prod_{k=2,4,8,\ldots,\frac n4}
\left(
{2}\cdot{2^{\vphantom{i}^{-\frac{n}{k}}}}
\KKNK
\right)^{\EExpLog[-0.5ex]{k}}
\label{LBru:eq:4}\\
&=&
                        \ExpHam{2}
\cdot
(   3^{\Eexp{4}}  \cdot \ExpHam{4} )
\nonumber\\&&
\qquad\qquad\,\,\,{}\cdot
(  35^{\Eexp{8}}  \cdot \ExpHam{8} )
\nonumber\\&&
\qquad\quad\,{}\cdot
(6435^{\Eexp{16}} \cdot \ExpHam{16} )
\nonumber\\&&
\,\cdot
\dots\cdot
\left(\left(
{\textstyle\frac 12{\big\lgroup\hspace{-0.7ex}}{{n/4}\atop{n/8}}
 {\hspace{-0.7ex}\big\rgroup}
}
\right)^{2^3}\cdot 2^{2^1}\right)
\cdot
\DmMinusOne
\nonumber
\end{eqnarray}
\end{theorem}
\begin{proof}
Число способов определить расширенный $1$-совершенный код по формулам (\ref{LBru:eq:2})
с ограничениями Конструкции~\ref{LBru:con:3}
легко посчитать при помощи Следствия~\ref{LBru:p:Dt} и Предложения~\ref{LBru:p:DtVt}.
Предложение~\ref{LBru:p:bold} гарантирует, что различные наборы локальных автоморфизмов дают различные коды.
\end{proof}
Поскольку имеется взаимнооднозначное соответствие (удаление последнего символа) между расширенными
 $1$-совершенными кодами и {\it $1$-совершенными кодами}, имеем следующее:

\begin{theorem}[Нижняя оценка]
Число $B(n-1)$
$1$-совершенных двоичных кодов длины $n-1=2^m-1$
удовлетворяет неравенству
\begin{equation}\label{LBru:eq:lbLA}
B(n-1)\geq \widetilde K_{\scriptscriptstyle LA}(n),
\end{equation}
где точное выражение и асимптотика для $\widetilde K_{\scriptscriptstyle LA}(n)$
представлены в Теореме~\ref{LBru:th:Kla}.
\end{theorem}

Как легко заметить, предыдущая оценка
$\ExpHam{2}\cdot(   3^{\Eexp{4}}  \cdot \ExpHam{4} )$
\cite{LBru:Kro:quas_l_b} состоит из двух множителей ($k=2,4$) из (\ref{LBru:eq:4}).


\begin{conjecture}\label{LBru:cj:} Оценка {\rm (\ref{LBru:eq:lbLA})} асимптотически точна,
т.\,е. {\rm (\ref{LBru:eq:4})} --- асимптотика числа $1$-совершенных двоичных кодов длины $n-1=2^m-1$.
\end{conjecture}

Эта гипотеза подтверждается имеющейся информацией об $1$-совершенных кодах малых рангов,
т.\,е. рангов $+1$ и $+2$. \emph{Ранг} (\emph{rank}) кода --- размерность его аффинной оболочки;
скажем, что $1$-совершенный код длины $n-1$ есть код ранга $+p$,
{{если}} его ранг равен $r_H+p$,
где $r_H$ --- размерность линейного $1$-совершенного кода (кода Хемминга) той же длины.
(Термин `аффинная оболочка' означает то же, что и `линейная оболочка'
для кодов, содержащих $\bar 0$, но аффинная оболочка инвариантна относительно сдвигов кода.)
Известно, что ЛА конструкция дает почти все коды ранга $+1$
и почти все коды ранга $+2$ (и, разумеется, коды большего ранга).
Более того, при фиксированной аффинной оболочке число
$1$-совершенных кодов ранга $+1$ асимптотически равно
$$ \ExpHam{2}, $$
а ранга $+2$ ---
$$
\ExpHam{2}\cdot
(   3^{\Eexp{4}}\cdot
    \ExpHam{4} )
    $$
(если не фиксировать аффинную оболочку кода, то эти значения нужно умножить на
$n!/2^{n/2}(\frac{n}{2}-1)(\frac{n}{2}-2)\ldots(\frac{n}{2}-\frac{n}{4})$ и
$n!/24^{n/4}(\frac{n}{4}-1)(\frac{n}{4}-2)\ldots(\frac{n}{4}-\frac{n}{8})$ соответственно).
Эти данные следуют из представления $1$-совершенных
двоичных кодов ранга $+1$ и $+2$ \cite{LBru:AvgHedSol:class}
и асимптотики $3^{n+1}2^{2^n +1}(1+o(1))$
числа $n$-арных квазигрупп порядка $4$ \cite{LBru:KroPot:Lyap},\cite{LBru:PotKro:asymp}.

\begin{remark}
Все коды из Конструкции~\ref{LBru:con:3} имеют дефицит ранга ${\rm RD}=2$
(максимальный возможный ранг расширенного $1$-совершенного кода длины $n\geq 16$ равен $n-1$, см. \cite{LBru:EV:94}; таким образом, дефицит ранга определяется как ${\rm RD}(C)\eqdf (n-1)-{\rm rank}(C)$)
и
(как следует из неравенства
$\dim({\rm kernel}(C))\geq 2^{{\rm RD}(C)}$
для двоичных $1$-совершенных кодов ранга $+2$ и больше, см. \cite[Corollary 2.6]{LBru:PheVil:RankKernel})
размерность ядра (kernel) не меньше $4$,
где ${\rm kernel}(C)\eqdf \{\bar k\,|\,C+\bar k=C\}$.
Последний факт означает, что конструкция дает не менее
$$ \frac{\widetilde K_{\scriptscriptstyle LA}(n)}{n!2^{n-5}}$$
неэквивалентных расширенных $1$-совершенных двоичных кодов длины $n$ и
$$ \frac{\widetilde K_{\scriptscriptstyle LA}(n)}{(n-1)!2^{n-5}}$$
неэквивалентных $1$-совершенных двоичных кодов длины $n-1$,
где $n!2^{n-1}$ --- число изометрий $\F$
и $(n-1)!2^{n-1}$ --- число изометрий $F^{n-1}$.
\end{remark}

Да, Гипотеза~\ref{LBru:cj:} подразумевает, что почти все (расширенные) $1$-совершенные коды имеют ранг $n-3$, который не является полным ($n-1$ для расширенных $1$-совершенных кодов длины $n$),
и даже не предполный ($n-2$).
Это не подтверждается для кодов длины $16$,
см. \cite{LBru:ZZ:2006:l16r14} (кстати, расширенных $1$-совершенных кодов длины $16$ полного ранга действительно мало, см. \cite{LBru:OstPot:16},\cite{LBru:OPP:16}),
но ЛА конструкция еще не ``набрала обороты'' при $n=16$ ($m=4$).
В самом деле, при $m=4$ среди трех множителей (\ref{LBru:eq:4}) первый ($t=1$)
имеет почти такую же величину, что и второй ($t=2$), а множитель $n!/6((n/4)!)^4$
является наибольшим, в то время как асимптотически первый множитель --- самый мощный.
С другой стороны, факт, что почти все коды не имеют полный ранг, вряд ли был бы сюрпризом.
Например, аналогичное явление известно для МДР-кодов с расстоянием $2$ в четырехбуквенном алфавите
($n$-арных квазигрупп порядка $4$, см. \cite{LBru:PotKro:asymp},\cite{LBru:KroPot:n-ary}), у которых
ранг (над $Z_2^2$) имеет три возможных значения
(ранг $n-1$ для линейных кодов длины $n$, ранг $n-\frac{1}{2}=\log_4|Z_2^{2n-1}|$ для `полулинейных',
и ранг $n$),
причем асимптотически почти все коды имеют средний ранг.
Стоит также заметить, что число неэквивалентных систем троек Штейнера порядка $16$
 полного ранга $15$ меньше, чем ранга $14$ \cite{LBru:KOP}.

\pagebreak
\label{LBru:sec:biblio}
\providecommand\href[2]{{#2}} \providecommand\url[1]{\href{#1}{#1}}